\documentclass[11pt]{amsart}   
\usepackage{graphicx, amsmath, amssymb, amsthm}
\usepackage{epsf}
\DeclareGraphicsExtensions{.jpg,.pdf,.mps,.png}
\newtheorem{question}{Question}
\theoremstyle{definition}

\theoremstyle{remark}

\numberwithin{equation}{section}

\begin{document}
\title{A note on fixed points of abelian actions in dimension one}

\author{J. P. Boro\'nski}
\address{National Supercomputing Centre IT4Innovations, Division of the University of Ostrava,
Institute for Research and Applications of Fuzzy Modeling,
30. dubna 22, 701 03 Ostrava,
Czech Republic -- and -- Faculty of Applied Mathematics,
AGH University of Science and Technology,
al. Mickiewicza 30,
30-059 Krak\'ow,
Poland}
\email{jan.boronski@osu.cz}



\begin{abstract}
The result of Boyce and Huneke gives rise to a 1-dimensional continuum, which is the intersection of a descending family of disks, that admits two commuting homeomorphisms without a common fixed point. 
\end{abstract}

\maketitle
Two commuting $C^1$-diffeomorphisms of a disk must have a common fixed point, and if they preserve a nonseparating plane continuum $Z$ (i.e. a compact and connected set such that $\mathbb{R}^2\setminus Z$ is connected), then there is a common fixed point in $Z$ \cite{Ribon}. An arc-like continuum is the one that is the inverse limit of arcs. Bing showed that every arc-like continuum embeds in the plane as the intersection of a descending family of disks \cite{Bing}. Hamilton showed that every arc-like continuum has the fixed point property \cite{Ha}. In particular, any $\mathbb{Z}$-action on such a continuum fixes a point. For more general abelian actions on arc-like continua no such result has been known, although it holds for example for dendrites and uniquely arcwise connected continua \cite{Shi}. 
For self-maps of the arc the examples of Boyce \cite{Boyce} and Huneke \cite{Huneke} provide two commuting surjective maps without a common fixed point\footnote{I am grateful to Benjamin Vejnar for bringing this result to my attention. An interested reader might want to consult \cite{Vejnar} for some related recent results.}. Noteworthy their result gives the following. 
\vspace{0.5cm}

\noindent
\textbf{Example.} There exists an arc-like continuum $X$ and a pair of homeomorphisms $F,G:X\to X$ such that $F\circ G=G\circ F$ and $\operatorname{Fix}(F)\cap \operatorname{Fix}(G)=\emptyset$. 
\vspace{0.3cm}

\noindent
The argument for the existence of the above example is quite short. However, we did not find such a result in the literature, and none of the experts we consulted had been aware of such a result before. After personal discussions during the 52nd Spring Topology and Dynamics Conference, held at Auburn University in March of 2018, we were encouraged to publish this note. 
\begin{proof}

	Let $f,g:[0,1]\to[0,1]$ be two commuting surjections without a common fixed point. Let $h=f\circ g$ and consider the inverse limit space 

	$$X=\lim_{\leftarrow}\{[0,1],h\}=\{(x_1,x_2,...)\in[0,1]^{\mathbb{N}}:h(x_{i+1})=x_i,\textrm{ for all }i\in\mathbb{N}\}.$$ Since $f$ and $g$ commute, all three maps induce maps on $X$, given by 

	$$H(x_1,x_2,...)=(h(x_1),h(x_2),....),$$ 

	$$F(x_1,x_2,...)=(f(x_1),f(x_2),....),$$

	and

	$$G(x_1,x_2,...)=(g(x_1),g(x_2),....).$$ 

	To see that $F$ and $G$ are homeomorphisms one can use the following, particularly elegant, argument from \cite{Mouron}. Since $H=F\circ G=G\circ F$, and $H$ is just the shift homeomorphism, the maps $F$ and $G$ are also homeomorphisms. Since $f$ and $g$ have no common fixed point, the same is true for $F$ and $G$. 

\end{proof}

We note that by Barge-Martin embedding theorem \cite{BM} both $F$ and $G$ extend (up to conjugacy) to homeomorphisms of the plane. Therefore the following question is of interest. 
\begin{question}
Can $F$ and $G$ be extended to two commuting planar homeomorphisms?
\end{question}
An affirmative answer to the above problem would prove the necessity of the $C^1$ setting in the generalization of the Cartwright-Littlewood fixed point theorem in \cite{Ribon}. We also note, that the example in \cite{Huneke} is constructed as a limit of piecewise linear functions, each of which has a well-defined derivative of constant absolute value greater than $3$ everywhere, outside of a finite set of critical points. This gives a good starting point to a potential modification, that would result in the continuum $X$ presented here, being homeomorphic to the pseudo-arc; see \cite{MT}.  
\begin{question}
Can $f$ and $g$ be modified so that the example holds on the pseudo-arc? 
\end{question}
The interested reader might want to consult \cite{Mouron} for related results.
\vspace{1cm}

\noindent
\textbf{Acknowledgments.} The author is supported by University of Ostrava grant lRP201824 "Complex topological structures" and the NPU II project LQ1602 IT4Innovations excellence in science.
\label{lastpage}
\bibliographystyle{amsplain}

\end{document}